\documentclass[leqno,12pt]{article}
\usepackage{amssymb,amsfonts}
\usepackage{amsmath,latexsym}

\usepackage{amstext,epic,eepic,epsf,pslatex}
\usepackage{graphicx}

\newcommand{\bepr}{{\em Proof} } 
\newcommand{\enpr}{\hfill \rule{.5em}{.5em}}

\newcommand{\R}{{\mathbb R}}

\newcommand{\Tr}{\hbox{Tr\,}}

\def\XXint#1#2#3{{\setbox0=\hbox{$#1{#2#3}{\int}$ }
\vcenter{\hbox{$#2#3$ }}\kern-.6\wd0}}

\newtheorem{defin}{Definition}[section] 
\newtheorem{prop}{Proposition}[section] 
\newtheorem{thm}{Theorem}[section] 
\newtheorem{lemma}{Lemma}[section]

\newtheorem{cor}{Corollary}[section]

\begin{document}

\title{Hard spheres dynamics: weak {\em vs} strong collisions}

\author{Denis Serre \\ \'Ecole Normale Sup\'erieure de Lyon\thanks{U.M.P.A., UMR CNRS--ENSL \# 5669. 46 all\'ee d'Italie, 69364 Lyon cedex 07. France. {\tt denis.serre@ens-lyon.fr}}}

\date{}

\maketitle

\begin{abstract}
We consider the motion of a finite though large number $N$ of hard spheres in the whole space $\R^n$. Particles move freely until they experience elastic collisions. We use our recent theory of Compensated Integrability in order to estimate how much the particles are deviated by collisions. Our result, which is expressed in terms of hodographs, tells us that only $O(N^2)$ collisions are significant. 
\end{abstract}

\section{The dynamical model of spheres}

We consider a set of $N$ identical spheres of radius $a>0$ moving in the whole space $\R^n$. The coordinate in the physical space $\R^{1+n}$ is denoted $x=(t,y)$ where $t$ is the time and $y\in\R^n$ the position. 
In practice $n=3$ and $N$ is of the size of the Avogadro number, but the analysis below is valid in every space dimension and for any cardinality. We may think of the spheres as particles $P_\alpha$ ($1\le\alpha\le N$) of mass $m>0$, though it is not essential for our purpose. Since our estimates below do not depend on this parameter, we normalize it to $m=1$.

The  velocity $v_\alpha\in\R^n$ of $P_\alpha$ remains constant between two consecutive collisions. In particular, the trajectory of a particle is a polygonal chain. 
A collision between two particles $P_{\alpha,\beta}$ occurs when their centers $y_{\alpha,\beta}$ approach to a distance $2a$~:
$$|y_\beta-y_\alpha|=2a.$$ 
In terms of the velocities $v_{\alpha,\beta}$ and $v_{\alpha,\beta}'$ before and after the collision, ``approach'' is expressed by the inequalities
\begin{equation}
\label{eq:incv}
(v_\beta-v_\alpha)\cdot(y_\beta-y_\alpha)<0,\qquad (v_\beta'-v_\alpha')\cdot(y_\beta-y_\alpha)>0.
\end{equation}
We shall make use of the following properties of collisions:
\begin{itemize}
\item The combined momentum is conserved:
\begin{equation}
\label{eq:cmom}
v_\alpha'+v_\beta'=v_\alpha+v_\beta.
\end{equation}
\item  The collisions are friction-less, meaning that the jump of velocity is orthogonal to the common tangent space to the particles:
\begin{equation}
\label{eq:normalimp}
v_\alpha'-v_\alpha=v_\beta-v_\beta'\parallel y_\beta-y_\alpha.
\end{equation}
\item The energy, which is only of kinetic nature, is conserved along the collision:
\begin{equation}\label{eq:cen}
|v_\alpha'|^2+|v_\beta'|^2=|v_\alpha|^2+|v_\beta|^2.
\end{equation}
\end{itemize}

Two important quantities emerge from the considerations above, namely the total mass $M=Nm=N$ and the total energy
$$E=\sum_\alpha\frac12\,|v_\alpha(t)|^2\equiv\sum_\alpha\frac12\,|v_\alpha(0)|^2,$$
which do not depend on the instant at which they are computed.  They are therefore determined by the initial data and we see them as {\em a priori} given. A third conserved quantity is the total momentum
$${\cal Q}=\sum_\alpha v_\alpha(t)\equiv \sum_\alpha v_\alpha(0),$$
which satisfies $|{\cal Q}|\le\sqrt{2ME}$\,. It allows us to define a mean velocity $w:=\frac{\cal Q}M\,$. Our results shall be expressed in terms of the root mean square and standard deviation of the velocity distribution
$$\bar v:=\sqrt{\frac{2E}M\,}\,,\qquad{\bf v}:=\sqrt{\bar v^2-|w|^2\,}\,.$$
These are constants of the motion.

\subsection{State of the art}

For generic initial data, hard spheres dynamics is well defined globally in time, and every collision involves exactly two particles; see Alexander's Master thesis \cite{Alex}, or Theorem 4.2.1 of \cite{CIP}.  The question of finiteness of the number of collisions was raised by Ya. Sinai \cite{Sin} and solved by Vaserstein \cite{Vas}, whose work was simplified by Illner \cite{Illn,Illde}. Their proofs argue by contradiction and therefore do not yield an explicit upper bound for the number of collisions. The only known bound was found recently by Burago \& al. \cite{BFK}, in the form
$$\#\{\rm collisions\}\le(32 N^{3/2})^{N^2}.$$
On the opposite side, Burdzy \& Duarte \cite{BD} exhibit an initial configuration of $N$ hard spheres for which the number of collisions in the whole history is larger than $\frac1{27}\,N^3$. This lower bound is soon improved by Burago \& Ivanov \cite{BuIv} in $2^{\lfloor\frac N2\rfloor}$ where $\lfloor\cdot\rfloor$ is the floor function. 

The number of collisions can thus be extremely large. 
Nevertheless, our Theorem \ref{th:main} provides a realistic bound $O(\epsilon^{-1}N^2)$ for those collisions that have a significant impact on the short-time dynamics, in the sense that $|v'-v|\ge\epsilon \bar v$ for a given threshold $\epsilon>0$. This solves a question raised in \cite{BuIv}:
\begin{quote}
{\em It seems that, if the number of collisions is ÒlargeÓ, then the overwhelming number of collisions are ÒinessentialÓ in the sense that they result in almost zero exchange of momenta, energy, and directions of velocities of balls. We will think about it tomorrow.}
\end{quote}

\subsection{Main result}

At a collision between two particles, each one experiences a jump of velocity $\delta v:=v'-v$~; the trajectory of the center of each sphere experiences a {\em kink}. Because of the conservation of momentum, the jumps of both particles compensate each other. We shall measure whether a collision is weak or significant, in terms of the ratio $|\delta v|/\bar v$. Our motivation is that at a macroscopic level, where a flow is described by thermodynamical variables, it is expected that the pressure be related with the number of significant collisions.

\bigskip

In what follows, a {\em universal constant} $c_n$ is a finite number which may depend upon the space dimension $n$, but does not upon the initial configuration. The same notation $c_n$ occurs in various places, but the constant may differ from one line to another.
\begin{thm}\label{th:main}
Consider $N$ hard spheres moving in the physical space $\R^n$. Let $\bar v$ be the root mean square velocity of the system. Let us  assume (generic) that the collision set is finite on every band $(0,\tau)\times\R^n$, and that the motion involves only binary collisions.

There exists a universal constant such that the following inequality holds true
\begin{equation}
\label{eq:estbin}
\sum_{\rm kinks}\left(\bar v\,|v'-v|+\,|v\wedge v'|\right)\le c_nN^2\bar v^2,
\end{equation}
where the sum runs over the particles and the collisions they experience.
\end{thm}

The inequality above tells us that if we neglect the {\em weak} collisions, for which either $|v'-v|<\!\!<\bar v$ or $|v\wedge v'|<\!\!<\bar v\,^2$, then in average: -- a given particle will experience $O(N)$ collisions, --  two given particles will collide together $O(1)$ many times. These big O involve only universal constants.

A slight improvement of (\ref{eq:estbin}) can be obtained by applying these estimates to the flow measured in the inertial frame in which the mean velocity (a constant of the motion) vanishes. In term of the standard deviation $\bf v$, we have the following Galilean-invariant estimate.
\begin{equation}
\label{eq:galinv}
\sum_{\rm kinks}|v'-v|\le c_nN^2{\bf v}.
\end{equation}

\paragraph{Interpretation.} Given a particle $P_\alpha$ (for $\alpha\in[\![1,N]\!]$), we may consider its hodograph, which is the curve parametrized by $t\mapsto v_\alpha(t)$. For cosmetic reasons, we prefer to consider the graph translated by $w$, denoted $H_\alpha$, which is the hodograph when we use the inertial frame in which the mean velocity vanishes. Because the map $t\mapsto v_\alpha(t)-w$ is piecewise constant,  $H_\alpha$ is actually defined as the polygonal chain passing through the consecutive values of $v_\alpha-w$. Each $H_\alpha$ has a length $\ell(H_\alpha)$. It sweeps also an area $A(H_\alpha)$ about the origin. Our estimate can be recast as
\begin{equation}
\label{eq:estlH}
\sum_\alpha\ell(H_\alpha)\le c_nN^2{\bf v},\qquad\sum_\alpha A(H_\alpha)\le c_nN^2{\bf v}\bar v.\end{equation}
In particular $\ell(H_\alpha)$ and $A(H_\alpha)$ are finite for every $\alpha$. The length of a typical hodograph is an $O(N{\bf v})$, and the area sweapt is an $O(N{\bf v}\bar v)$. This gives us an information about how much the particles paths can become random as $N\rightarrow+\infty$.

\bigskip

\paragraph{Comments.}
\begin{itemize}
\item These estimates are independent of the radius $a$ of the particles. This was predictible, since a motion of particles of size $a$ yields a similar motion with particle size $a'$ by a change of scale. The number of collisions remains the same.
\item The estimates are also independent of the space  dimension, apart for the constants $c_n$.
\item In a previous version of this paper, we presented the weaker estimate 
$$\sum_{\rm kinks}\frac{\bar v^2\,|v'-v|^2+\,|v\wedge v'|^2}{\sqrt{(\bar v^2+\,|v|^2)(\bar v^2+\,|v'|^2)}}\le c_nN^2\bar v^2.$$
The improvement in (\ref{eq:estbin}) relies upon an explicit calculation of {\em determinantal masses} (see Section \ref{ss:Dm}), which replaces a coarse lower bound.
\item The quantity $|v\wedge v'|$ can be rewritten $|v|\cdot|v'|\sin\theta$ where $\theta$ is the angle by which the trajectory is deviated.
\item Inequality (\ref{eq:estlH}.1) is sharp in the sense that there exists an initial configuration for which ${\bf v}=1$ and $\sum_\alpha \ell(H_\alpha)=N^2$. Just consider $N=2p$ particles along a line, with initial velocities $\pm1$, $p$ particles moving to the left and the other $p$ to the right.
\item Since $\ell(H_\alpha)$ is finite, every particle $P_\alpha$ admits limit velocities $v_{\alpha\pm}$ as $t\rightarrow\pm\infty$, and we have the estimate
\begin{equation}
\label{eq:scat}
\sum_{\alpha=1}^N|v_{\alpha\pm}-v_\alpha(0)|\le c_nN^2{\bf v}.
\end{equation}
Contrary to (\ref{eq:estlH}.1), it is unlikely that (\ref{eq:scat}) be sharp.
\end{itemize}

\bigskip

\paragraph{Open questions.} 1) A more natural context arises when the spheres evolve in a bounded  domain $\Omega$ with impermeable wall. As far as we know, the finiteness of the set of collisions on finite time intervals is still an open question. At least, it is known that  a single particle can bounce infinitely many times at the boundary $\partial\Omega$ in finite time. A natural question is therefore to estimate the number and strengths of the collisions in a given compact subdomain $K$. The best that we could expect is that 
$$\frac1{|J|}\,\sum_{{\rm in}\,J\times K} |v'-v|$$
be uniformly bounded in terms of ${\rm dist}(K,\partial\Omega)$, $N$ and $\bar v$, when the length of the time interval $J$ is larger than the characteristic time ${\rm diam}\,\Omega/\bar v$. So far, we did not succeed to adapt our method, or to establish such a bound.
The analogous question arises for a space-periodic flow and remains open as well. 2) (Courtesy of L. Saint-Raymond) The construction of the mass-momentum tensor (Section \ref{s:MMT}) depends heavily on our assumption that the particles are spheres. What can be said of the dynamics of more general solids ? One might  consider  identical rigid convex bodies.

\paragraph{Plan of the paper.}
The central object of this paper is the {\em mass-momentum tensor} associated with the motion. Its construction, done in Section \ref{s:MMT}, is  more involved than that for flows obeying the Euler or the Boltzmann equations (see \cite{Ser_IHP}). Nevertheless the idea is the same, the entries representing the distribution of mass, momentum and stress. The tensor is a positive semi-definite symmetric matrix of size $1+n$, whose entries are bounded measures supported by a graph. The conservation of mass and momentum is expressed by the row-wise identity ${\rm Div}\,T=0$. One striking feature in this construction is the introduction of massless virtual particles ({\em collitons}) whose role is to carry the interchange of momentum between colliding spheres. The short Section \ref{s:CI} recalls the principles of Compensated Integrability for divergence-controlled positive tensors, as developped in our former papers \cite{Ser_IHP,Ser_CI}. Section \ref{s:Tgraph} is an improvement of the theory when such tensors are supported by graphs~; this is where we introduce the concept of determinantal mass at the nodes. With this tool in hands, the proof of Theorem \ref{th:main} becomes rather short and is carried out in
Section \ref{s:bin}. We apply the extended version of Compensated Integrability to a combination of the mass-momentum tensor and an appropriate parametrized complement.

\paragraph{Acknowledgement.} I am indebted to Laure Saint-Raymond and Reinhard Illner for valuable discussions and their help in gathering the relevant literature. \'Etienne Ghys remarked that my results can be rephrased in terms of hodographs.

\section{The mass-momentum tensor}\label{s:MMT}

From now on, we denote $d=1+n$ the time-space dimension. If $J$ is a line or a segment, we denote $\delta_J$ the one-dimensional Lebesgue measure along $J$. We recall that for a distribution $f$, positive homogeneity of a given degree $\kappa$ can be defined either by duality, or by the Euler identity $(x\cdot\nabla)f=\kappa f$. The Lebesgue measure over a $k$-dimensional linear subspace of $\R^d$ is homogeneous of degree $k-d$~; for instance ${\cal L}^d$ has degree $0$ (obvious), while $\delta_0$ has degree $-d$. If $L$ is a line, or a semi-line from the origin, its one-dimensional Lebesgue measure $\delta_L$ has degree $1-d$.

If $S$ is a $d\times d$ symmetric tensor over an open domain of $\R^n$, whose entries are distributions, the row-wise divergence ${\rm Div}\,S$ is the vector field whose entries are the distributions
$$\sum_{\beta=1}^d\partial_\beta S_{\alpha\beta},\qquad\alpha=1,\ldots,d.$$

If $Q\in \R^d\setminus\{0\}$, and $\eta=\frac Q{|Q|}$\,, then for every line $L=\bar x+\R\eta$, we form the symmetric tensor
$$S^{Q,L}:=Q\otimes \eta\,\delta_L.$$
In other words
$$\langle S_{ab},\phi\rangle=Q_a\eta_b\int_\R\phi(\bar x+s\eta)\,ds,\qquad\forall\phi\in C_K(\R^{1+n}),$$
where $\bar x$ is an arbitrary point on the line $L$. When the context makes it clear, we write $S^Q$ instead.
\begin{lemma}
\label{l:SQ}
One has
$${\rm Div}\,S^{Q,L}=0.$$
 \end{lemma}
 
 \bigskip
 
 \bepr
 
 If $\phi$ is a test function, then
 $$\langle{\rm Div}\,S^{Q,L},\phi\rangle=-\langle S^{Q,L},\nabla\phi\rangle=-Q\int_\R \eta\cdot\nabla\phi(\bar x+s\eta)\,ds=-Q\int_\R \frac d{ds}\,\phi(\bar x+s\eta)\,ds=0.$$
 
 \enpr

\subsection{Single particle}

We begin by considering a single particle $P$ of unit mass, whose constant velocity is $v\in\R^n$. The trajectory $t\mapsto (t,y(t))$ of the center of mass in the physical space $\R^{1+n}$ is a line $L$, whose direction is
$$\xi=\frac V{|V|}\,,\qquad\hbox{where } V:=\binom1v.$$
We define the {\em mass-momentum} tensor of $P$ to be $S^V$. According to Lemma \ref{l:SQ}, it is divergence-free.

\subsection{Multi-line configuration}

When $L$ is replaced by a semi-infinite line $L^+=\bar x+\R_+\eta$, the tensor
$$S^{Q+}:=Q\otimes \eta\,\delta_{L^+}$$
is no longer divergence-free. The calculation above yields
\begin{equation}
\label{eq:HL}
\langle{\rm Div}\,S^{Q+},\phi\rangle=\phi(\bar x)Q,
\end{equation}
which is recast in distributional terms as
$${\rm Div}\,S^{Q+}=Q\,\delta_{\bar x}\,.$$

Now, if finitely many vectors $Q_1,Q_2,\ldots$ are given, together with a point $\bar x\in\R^d$, we may form the converging semi-lines $L_j^+=\bar x +\R_+Q_j$ and define a symmetric tensor
$$S_{\rm multi}:=\sum_jS^{Q_j+}.$$
Then (\ref{eq:HL}) tells us that $S_{\rm multi}$ is divergence-free whenever
\begin{equation}
\label{eq:trpl}
\sum_jQ_j=0.
\end{equation}

\paragraph{Application to $1$-D dynamics.} When $n=1$, we may simplify the hard sphere model by setting $a=0$. 
At a binary collision, the particles meet at some point $\bar x\in\R^{1+1}$, with incoming velocities $v,w$ and outgoing ones $v',w'$. Let us choose 
\begin{equation}
\label{eq:VVVV}
V_1=-\binom1v,\qquad V_2=-\binom1{w},\qquad V_3=\binom1{v'},\qquad V_4=\binom1{w'}.
\end{equation}
Then the positive semi-definite tensor
$$T:=S^{V_1+}+S^{V_2+}+S^{V_3+}+S^{V_4+}$$
associated with this pair of particles is divergence-free~; the compatibility condition (\ref{eq:trpl}) is ensured by the conservation of mass and momentum through the collision. The support of $T$ is the union of the trajectories.

\subsection{Binary collisions ($n\ge2$)}\label{ss:bin}

When $n\ge2$ instead, the radius $a$ must be positive, in order that collisions take place. 

Let two particles $P_i$ and $P_j$ collide at some time $t^*$. The trajectory of $P_i$ displays a kink at a point $\bar x_i=(t^*,\bar y_i)$, and that of $P_j$ does at $\bar x_j=(t^*,\bar y_j)$ at the same instant $t^*$. We have $|\bar y_j-\bar y_i|=2a$. Let us define $V_1,\ldots,V_4$ as in (\ref{eq:VVVV}). Locally, the trajectories are made of segments of  the semi-lines
$$L_1^+=\bar x_i+\R_+V_1,\qquad L_2^+=\bar x_j+\R_+V_2,\qquad L_3^+=\bar x_i+\R_+V_3,\qquad L_4^+=\bar x_j+\R_+V_4.$$
Because the lines do not meet at a single point, the tensor $S=S^{V_1+}+S^{V_2+}+S^{V_3+}+S^{V_4+}$ is not divergence-free. We have instead
$${\rm Div}\,S=(V_1+V_3)\delta_{x_i}+(V_2+V_4)\delta_{x_j}=(V_1+V_3)(\delta_{x_i}-\delta_{x_j}).$$

In order to recover a divergence-free tensor, we introduce the vector $Q$
$$Q=\binom0q,\qquad q=v'-v=w-w'.$$
Because of (\ref{eq:normalimp}), the segment $C=[\bar x_i,\bar x_j]$ has direction $Q$. In the neighbourhood of the collision, we can define the tensor 
$$T=S^{V_1+}+S^{V_2+}+S^{V_3+}+S^{V_4+}+S^{Q,C}.$$
Each of the five terms in the sum above is divergence-free away from either $\bar x_i$ or $\bar x_j$. At $\bar x_i$, ${\rm Div}\,T$ is a sum of three Dirac masses, whose weight is
\begin{equation}
\label{eq:VVQ}
V_1+V_3-Q=\binom{-1+1-0}{-v+v'-q}=0,
\end{equation}
where the minus sign in front of $Q$ comes from the fact that $Q$ is oriented from $x_j$ to $x_i$.
A similar identity holds true at $\bar x_j$, with now a plus sign in front of $Q$. We conclude that 
$${\rm Div}\,T=0.$$

We may interpret the contribution $S^{Q,C}=Q\otimes\eta\,\delta_C$ as that of a virtual particle. This particle is  massless, because the first component of $Q$ vanishes. It carries the momentum which is exchanged instantaneously between $P_i$ and $P_j$. We suggest the name {\em colliton} for this object.

\subsection{The complete construction}

Assuming again that only binary collisions occur, we consider the union of trajectories of the centers of the $N$ particles. Each trajectory is a polygonal chain whose kinks occur where and when the particle suffers a collision. We {\bf define} the {\em mass-momentum} tensor $T$ of the configuration as the sum of the following contributions:
\begin{itemize}
\item For each segment $J$ of a trajectory between two consecutive collisions, the tensor
$$S^{V,J}=V\otimes\xi\,\delta_J,\qquad V=\binom1v,\qquad\xi=\frac V{|V|}\,,$$
where $v$ is the particle velocity along $J$.
\item For each binary collision, the corresponding colliton, as described in the previous paragraph.
\end{itemize}

The mass-momentum tensor $T$  is a divergence-free symmetric positive semi-definite tensor. We point out that its support is a graph, a one-dimensional object in $\R^{1+n}$. Thus $T$ vanishes almost everywhere in the Lebesgue sense. The support can be equiped with the positive measure $\Tr\,T$, with respect to which $T$ is rank-$1$ almost everywhere.

\bigskip

\paragraph{Finiteness.}
Because the particles are finitely many, and the collisions are finitely many in every band $H_\tau=(0,\tau)\times\R^n$ by assumption, the restriction of the entries of $T$ to $H_\tau$ are finite measures. This property is an essential hypothesis in Compensated Integrability (Theorem \ref{th:CI} below). Remark however that we do not have a practical bound of the total mass of $T$ in $H_\tau$, because we do not control efficiently the number and the strength of the collitons.

\section{Compensated integrability}\label{s:CI}

We shall make use of our recent theory of Compensated Integrability for divergence-controlled positive symmetric tensors, for which we refer to \cite{Ser_IHP,Ser_CI}. The appropriate version is given in the theorem below. Let $U\subset\R^d$ be an open set. Let $S$ be a distribution over $U$ that takes values in the cone ${\bf Sym}_d^+$. By positiveness, the entries $S_{ab}$ are locally finite measures. We say that $S$ is {\em
divergence-controlled} if these entries, as well as the coordinates of ${\rm Div}\,S$, are  finite measures. We recall that a divergence-controlled tensor admits a normal trace $S\vec\nu$ along the boundary $\partial U$, which {\em a priori} belongs to the dual space of ${\rm Lip}(\partial U)$.

We denote $\|\mu\|$ for the total mass of a (vector-valued) bounded measure $\mu$,
$$\|\mu\|=\langle|\mu|,{\bf1}\rangle.$$
This notation is used below in two distinct contexts, whether $\mu$ is a measure over a $(1+n)$-dimensional slab $H=(t_-,t_+)\times\R^n$, or a measure over $\R^n$. 
\begin{thm}\label{th:CI}
Let $H=(t_-,t_+)\times\R^n$ be a slab in $\R\times\R^n$, and $S$ be symmetric positive semi-definite tensor defined over $H$. We assume that $S$ is divergence-controlled\footnote{Mind that ${\rm Div}\,S$ involves time and space derivatives.} in $H$. Finally we assume that
the normal traces $S\vec e_t$ at the initial and final times $t=t_\pm$ are themselves bounded measures.

Then the measure $(\det S)^{\frac1{n+1}}$ actually belongs to $L^{1+\frac1n}(H)$ and we have
\begin{equation}
\label{eq:CIn}
\int_H(\det S)^{\frac1n}dy\,dt\le c_n\left(\|S\vec e_t(t_-)\|+\|S\vec e_t(t_+)\|+\|{\rm Div}\,S\|\right)^{1+\frac1n},
\end{equation}
where $c_n$ is a finite constant independent of $S$ and $H$.
\end{thm}

\paragraph{Remarks.} 
\begin{itemize}
\item The additional assumption that the normal traces are  bounded measures is equivalent to saying that  the extension $\tilde S$ by $0_{1+n}$ away from $H$ enjoys  the property that ${\rm Div}\,\tilde S$ is a bounded measure. We then have the formula
$${\rm Div}\,\tilde S=\widetilde{{\rm Div}\,S}-S\vec e_t(t_-)\otimes\delta_{t=t_-}+S\vec e_t(t_+)\otimes\delta_{t=t_+}\,,$$
where the first term in the right-hand side is again the extension by zero of ${\rm Div}\,S$ away from $H$.
\item
This theorem is useless when $S$ is rank-$1$ almost everywhere, because then $(\det S)^{\frac1d}\equiv0$ and the estimate (\ref{eq:CIn}) is trivial. The goal of the next section is to improve the statement when the tensor is supported by a graph.
\end{itemize}

\section{Tensors supported by a graph}\label{s:Tgraph}

Let $G$ be a non-oriented graph included in $H$, with straight edges. Let $S$ be a tensor of the form
\begin{equation}\label{eq:Sgraph}
S=\sum_Ja_J\eta_J\otimes\eta_J\delta_J,
\end{equation}
where the sum runs over the edges. The unit vector $\pm\eta_J$ is the direction of $J$, and $a_J>0$ is a weight on the edge. We already know that 
$${\rm Div}\,S=\sum_wm(w)\delta_w$$
where the sum runs over the vertices and the weight is given by 
$$m(w)=\sum_{J\sim w} a_J\eta_J.$$
This sum runs over the edges around the vertex $w$, with $\eta_J$ oriented outward. 

The tensor $S$ is positive semi-definite, its entries being locally finite measures. The divergence is a finite measure too, as well as the normal traces at $t=t_\pm$. For instance
$$S\vec e_t(t_-)=\sum' a_J\eta_J\delta_{x_{J-}}$$
where the sum runs over the set of edges that meet the hyperplane $t=t_-$, and $(t_-,x_{J-})$ is the intersection of this space with $J$.

When applying Theorem \ref{th:CI}, we have therefore a good control of the right-hand side. But as mentionned above, the estimate is useless because the left-hand side vanishes identically, due to $(\det S)^{\frac1d}\equiv0$. We shall see below that something can be gained at those vertices where $m(w)=0$ ($S$ is locally divergence-free), provided that the set of directions $\{\eta_J\,:\,J\sim w\}$ span $\R^d$.

\subsection{Minkowski potentials}

To begin with, we recall that if $U\subset\R^d$ is a convex open subset and $\theta\in W^{2,d-1}(U)$ is given, then the cofactor matrix $\Lambda^\theta:=\widehat{{\rm D}^2\theta}$ of the Hessian is symmetric, integrable and divergence-free. If in addition $\theta$ is convex, then this tensor is positive semi-definite. Because of the formula 
$$\det\widehat R=(\det R)^{d-1}$$
for $d\times d$ matrices, the expression 
$$\int_U(\det \Lambda^\theta)^{\frac1{d-1}}dx$$
which is at stake in Compensated Integrability equals
$$\int_U\det {\rm D}^2\theta\,dx={\rm vol}_d(\nabla\theta(U)).$$

The Sobolev regularity of the potential can actually be lowered, and one shows that every convex $\theta$ yields a non-negative divergence-free tensor $\Lambda^\theta$. Of special interest is the case where $\theta$ is positively homogeneous of degree $1$ (for instance, $\theta$ might be a norm)~; see \cite{Ser_CI}. Then $\Lambda^\theta$ is positively homogeneous of degree $1-d$,
$$\Lambda^\theta=\mu^\theta\left(\frac x{|x|}\right)\,\frac{x\otimes x}{|x|^{d+1}}\,,$$
where $\mu^\theta$ is some positive finite measure over the unit sphere $S_{d-1}$. This measure satisfies the relation
$$\int_{S_{d-1}}e\,\mu^\theta(e)=0.$$
When $\theta\in W^{2,d-1}_{\rm loc}$, $\mu^\theta$ is just an integrable function.

\bigskip

Conversely, if $\mu$ is a positive measure over $S_{d-1}$, one may form the positive symmetric tensor
\begin{equation}\label{eq:Mmu}
\Lambda=\mu\left(\frac x{|x|}\right)\,\frac{x\otimes x}{|x|^{d+1}}\,.
\end{equation} 
In other words, for every test function $\phi\in{\cal C}_K(\R^d))$, we have
$$\langle\Lambda_{ab},\phi\rangle=\int_0^\infty dr\int_{S_{d-1}}\omega_a\omega_b\phi(r\omega)\mu(\omega).$$
The tensor $\Lambda$ turns out to be divergence-free if and only if
\begin{equation}\label{eq:mue}
\int_{S_{d-1}}e\,\mu(e)=0.
\end{equation} 
The problem of whether there exists a convex potential $\theta$, positively homogeneous of degree $1$, such that $\Lambda=\Lambda^\theta$ received a positive answer, given by Pogorelov \cite{Pog}. The solution exists and is unique, up to the addition of a linear form. When the support of $\mu$ spans $\R^d$, this problem is equivalent to that of Minkowski, which asks for a convex body whose Gaussian curvature (here $\mu$) is prescribed as a function of the unit normal. For this reason, we call $\theta$ the {\em Minkowski potential} of $\Lambda$ (or of $\mu$).

The special case where $\mu$ is a finite sum of Dirac masses is precisely that solved by Minkowski himself \cite{Min} ; then the body is a convex polytope. We shall use only this case below.

\subsection{Determinantal mass at a vertex}\label{ss:Dm}

\begin{defin}
Let $\theta:\R^d\rightarrow\R$ be a convex function, positively homogeneous of degree one. Then the {\em determinantal mass} of $\Lambda^\theta$ at the origin is the volume of the convex body enclosed by the boundary $\nabla\theta(S_{d-1})$. It is denoted ${\rm Dm}(\Lambda^\theta;0)$.

If a divergence-controlled positive semi-definite tensor $S$ coincides with $\Lambda^\theta$ in a neighbourhood of the origin, we define again ${\rm Dm}(S;0):={\rm Dm}(\Lambda^\theta;0)$. Finally, using translations, we define the determinantal mass ${\rm Dm}(S;x^*)$ of a tensor $S$ at an arbitrary point $x^*\in\R^d$, provided $S$ coincides locally with $\Lambda^\theta(\cdot-x^*)$.
\end{defin}

The definition above applies in particular to the following situation. Let $S$ be a divergence-controlled positive semi-definite symmetric tensor over  a band $H$, and suppose that in a neighbourhood $U$ of some point $x^*$, it is of the form (\ref{eq:Sgraph}) for finitely many edges attached to the vertex $x^*$. Assume also that $S$ is divergence-free in this neighbourhood, that is $m(x^*)=0$.

Up to a translation, we may assume $x^*=0$. Then $S$ is locally homogeneous of degree $1-d$, and the equation $m(0)=0$ just says that $S$ fulfills condition (\ref{eq:mue}). By Pogorelov's theorem, it can therefore be parametrized locally as $S=\Lambda^\theta$ for some convex function, positively homogeneous of degree $1$. In particular $\theta$ is not smooth at the origin.

We use the determinantal mass as follows. Let us smooth out $\theta$ in a ball $B$ such that $2B\subset U$. The resulting convex potential $\xi$ is $C^\infty$ in $B$ and coincides with $\theta$ in the corona $2B\setminus B$. The associated tensor $\Lambda^\xi$ is $C^\infty$ in $B$ and coincides with $\Lambda^\theta$ in $2B\setminus B$. Thus we may form the divergence-controlled tensor $\tilde S$ such that $\tilde S=S$ in $H\setminus B$ and $\tilde S=\Lambda^\xi$ in $2B$. We have ${\rm Div}\,\tilde S={\rm Div}\,S$ because on the one hand $\tilde S=S$ away from $B$, and on the other hand both of them vanish in $2B$. Besides, the normal traces at $t=t_\pm$ coincide. When applying (\ref{eq:CIn}) to $\tilde S$, the right-hand side is therefore unchanged. However, the left-hand side gains the contribution
$$\int_B(\det \Lambda^\xi)^{\frac1{d-1}}dy\,dt=\int_B\det{\rm D}^2\xi\,dy\,dt={\rm vol}(\nabla\xi(B)).$$
Because $\nabla\xi\equiv\nabla\theta$ over $\partial B$, the latter quantity is precisely the determinantal mass of $S$ at the vertex. Applying Estimate (\ref{eq:CIn}) to $\tilde S$, we obtain
\begin{eqnarray*}
\int_H(\det S)^{\frac1n}dy\,dt+{\rm Dm}(S;x^*) & = & \int_H(\det \tilde S)^{\frac1n}dy\,dt \\
& \le & c_n\left(\|\tilde S\vec e_t(t_-)\|+\|\tilde S\vec e_t(t_+)\|+\|{\rm Div}\,\tilde S\|\right)^{1+\frac1n} \\
&  & = c_n\left(\|S\vec e_t(t_-)\|+\|S\vec e_t(t_+)\|+\|{\rm Div}\,S\|\right)^{1+\frac1n}.
\end{eqnarray*}

More generally, applying the construction described above at every vertex where $S$ is graph-like and divergence-free, we obtain the following improvement of Theorem \ref{th:CI}.
\begin{thm}\label{th:vrtx}
Let $H=(t_-,t_+)\times\R^n$ be a slab in $\R\times\R^n$, and $S$ be a symmetric, positive semi-definite tensor defined over $H$. We assume that $S$ is divergence-controlled and that
the normal traces $S\vec e_t$ at the initial and final times $t=t_\pm$ are bounded measures too.

Then we have
\begin{equation}
\label{eq:CIvrtx}
\int_H(\det S)^{\frac1n}dy\,dt+\sum{\rm Dm}(S;x^*)\le c_n\left(\|S\vec e_t(t_-)\|+\|S\vec e_t(t_+)\|+\|{\rm Div}\,S\|\right)^{1+\frac1n},
\end{equation}
where the summation extends over the vertices $x^*\in H$ about which $S$ is of the form (\ref{eq:Sgraph}) and is divergence-free.
\end{thm}

\paragraph{Remarks.}
\begin{itemize}
\item The calculation above suggests to redefine $(\det S)^{\frac1{d-1}}$ as the sum of an absolutely continuous part, the one at stake in Theorem \ref{th:CI}, and a singular one, made of Dirac masses ${\rm Dm}(S;x^*)\,\delta_{x^*}$ at every point $x^*$ where $S$ is graph-like and divergence-free.
\item The map $S\mapsto{\rm Dm}(S;x^*)$ is homogeneous of degree $\frac d{d-1}$ \,, invariant under the action of the orthogonal group: if $R$ is a rotation and $S$ is given by (\ref{eq:Sgraph}), then ${\rm Dm}(S;x^*)={\rm Dm}(S^R;x^*)$, where $S^R$ is defined by rotating each of the $\eta_J$'s and keeping $x^*$ and $a_J$ unchanged.
\item There is nothing special in the choice of a slab. Theorem \ref{th:vrtx} has a version in an arbitrary bounded open domain $\Omega\subset\R^d$.
\end{itemize}

\subsection{Calculations of determinantal masses}

When $\Lambda$ is given as in (\ref{eq:Mmu}), we do not know of a closed form of ${\rm Dm}(\Lambda;0)$ in terms of the measure $\mu$. We shall use only simple cases where the Minkowski potential and the determinantal mass can be calculated explicitly.
We present below two useful situations.

\subsubsection{The planar case}

When $d=2$ (that is $n=1$), our closed formula uses the linearity of the operator  $\theta\mapsto\widehat{{\rm D}^2\theta}$.
\begin{prop}[$d=2$]\label{p:closedF}
Let $\mu$ be a positive measure over the unit circle $S_1$ satisfying the constraint (\ref{eq:mue}). Let $\Lambda$ be the divergence-free tensor defined over $\R^2$ by (\ref{eq:Mmu}). Then the Minkowski potential of $\Lambda$ is given by
$$\theta(x)=|x|\,p\left(\frac x{|x|}\right),$$
where $p$ is a $2\pi$-periodic solution of $p+p''=\mu$ (the derivatives are taken with respect to angle).

The determinantal mass of $\Lambda$ at the origin is given by
\begin{equation}
\label{eq:FVWZ}
{\rm Dm}(\Lambda;0)=\frac18\int_0^{2\pi}\!\int_0^{2\pi}\mu(s_1)\mu(s_2)\sin|s_2-s_1|\,ds_1ds_2.
\end{equation}
\end{prop}

We point out that (\ref{eq:mue}), which writes here
$$\int_0^{2\pi}\mu(s)\sin s\,ds=\int_0^{2\pi}\mu(s)\cos s\,ds=0,$$
is precisely the solvability condition of $p+p''=\mu$ in the realm of periodic functions.

\bigskip

We shall use the following consequence of Proposition \ref{p:closedF}.
\begin{cor}\label{c:VWZ}
Let $V,W,Z\in\R^2$ be given, such that $V+W+Z=0$. For some $x^*\in\R^2$, and three converging lines of directions $V,W$ and $Z$, we consider the divergence-free tensor
$$\Lambda:=S^V+S^W+S^Z.$$
We have
\begin{equation}
\label{eq:DmVWZ}
{\rm Dm}(\Lambda;x^*)=\frac14\,|\det(V,W)|.
\end{equation}
\end{cor}

\bigskip

Mind that this expression is symmetric in $V,W,Z$. For instance,
$$\det(V,Z)=\det(V,-V-W)=-\det(V,W).$$
Actually, if we label $V,W,Z$ in the trigonometric order, with arguments $0\le\alpha<\beta<\gamma<2\pi$, the measure $\mu$ is given by
$$\mu=|V|\delta_\alpha+|W|\delta_\beta+|Z|\delta_\gamma,$$
and the formula (\ref{eq:FVWZ}) gives
\begin{eqnarray*}
{\rm Dm}(\Lambda;x^*) & = & \frac14(|V|\cdot|W|\,\sin(\beta-\alpha)+|W|\cdot|Z|\,\sin(\gamma-\beta)+|Z|\cdot|V|\,\sin(\gamma-\alpha)) \\
& = & \frac14\,(|\det(V,W)|+|\det(W,Z)|-|\det(Z,V)|)=\frac14\,|\det(V,W)|,
\end{eqnarray*}
where we have used $\beta-\alpha,\,\gamma-\beta\in(0,\pi)$ and $\gamma-\alpha\in(\pi,2\pi)$.

\bigskip

\bepr 
(of Proposition \ref{p:closedF}.)

 Let $\theta$ be the Minkowski potential of $\Lambda$ and $p$ be its restriction to the unit circle. An elementary calculation yields
$$\Lambda=\widehat{{\rm D}^2\theta}=\begin{pmatrix} \theta_{,22} & - \theta_{,12} \\ - \theta_{,12} &  \theta_{,11} \end{pmatrix}=(p+p'')\,\frac{x\otimes x}{|x|^3}\,,$$
whence the differential equation $p+p''=\mu$. As mentionned above, there exists a $2\pi$-periodic solution $p$ because of the constraints (\ref{eq:mue}) which express the divergence-freeness of $\Lambda$. The solution is unique up to the addition of $a\sin+b\cos$~; in terms of $\theta$, this means uniqueness up to the addition of a linear form. The image of $\nabla\theta$ is the curve
$$\phi\longmapsto\binom{p(\phi)\cos\phi-p'(\phi)\sin\phi}{p(\phi)\sin\phi+p'(\phi)\cos\phi}.$$
The mass ${\rm Dm}(\Lambda;0)$, being the area enclosed by this curve, equals
\begin{eqnarray*}
{\rm Dm}(\Lambda;0) & = & \frac12\,\int_0^{2\pi}\theta_{,1}d\theta_{,2}=\frac12\int_0^{2\pi}(p\cos\phi-p'\sin\phi)(p+p'')\cos\phi\,d\phi \\
& = & \frac12\int_0^{2\pi}(p\cos\phi-p'\sin\phi)\mu\cos\phi\,d\phi.
\end{eqnarray*}
Let us define
$$\lambda(\phi):=\int_0^\phi\mu(s)\cos s\,ds,$$
which is periodic because of (\ref{eq:mue}). Then
\begin{eqnarray*}
{\rm Dm}(\Lambda;0) & = & \frac12\int_0^{2\pi}(p\cos\phi-p'\sin\phi)\lambda'\,d\phi =\frac12\int_0^{2\pi}(p+p'')\lambda\sin\phi\,d\phi=\frac12\int_0^{2\pi}\mu\lambda\sin\phi\,d\phi \\
& = & \frac12\int_0^{2\pi}d\phi\int_0^\phi\mu(\phi)\mu(s)\sin\phi\cos s\,ds.
\end{eqnarray*}
Using again (\ref{eq:mue}), this gives
$${\rm Dm}(\Lambda;0)=-\frac12\int_0^{2\pi}d\phi\int_\phi^{2\pi}\mu(\phi)\mu(s)\sin\phi\cos s\,ds.$$
With Fubini, this yields
$${\rm Dm}(\Lambda;0)=-\frac12\int_0^{2\pi}\int_0^s\mu(\phi)\mu(s)\sin\phi\cos s\,d\phi\,ds\stackrel{\rm relabel.}{=\!=}-\frac12\int_0^{2\pi}\int_0^\phi\mu(\phi)\mu(s)\sin s\cos\phi\,ds\,d\phi.$$
Combining the formul\ae\, above, we obtain
$${\rm Dm}(\Lambda;0)=\frac14\int_0^{2\pi}d\phi\int_0^\phi\mu(\phi)\mu(s)\sin(\phi- s)\,ds,$$
which by symmetrization, yields (\ref{eq:FVWZ}).

\enpr

\bigskip

\subsubsection{Direct sums}

We continue our study of tensors of the form (\ref{eq:Sgraph}), say centered at the origin. We suppose here that the set of vectors $\eta_J$ can be split into two subsets, orthogonal to each other: $\R^d=E_-\oplus^\bot E_+$, and each $\eta_J$ is either in $E_-$ or in $E_+=E_-^\bot$. Up to a rotation, we may always assume that $E_-=\R^p\times\{0\}$ and $E_+=\{0\}\times\R^q$ with $p+q=d$. Our tensor writes therefore blockwise 
\begin{equation}\label{eq:Spm}
S=\begin{pmatrix} S_-\otimes\delta_{x_+=0} & 0 \\ 0 & \delta_{x_-=0}\otimes S_+ \end{pmatrix}.
\end{equation}
The tensors $S_\pm$ are defined over open subsets of $E_\pm$ and inherit the divergence-freeness. They are actually of a form similar to (\ref{eq:Sgraph}), though in either $\R^p$ or $\R^q$ instead of $\R^d$. Each of both admits a Minkowski potential:
$$
S_-=\widehat{{\rm D}_-^2\theta_-},\qquad S_+=\widehat{{\rm D}_+^2\theta_+}
$$
where $\theta_\pm$ is a convex function of $x_\pm$, positively homogeneous of degree $1$. The derivative $D_-$ (resp. $D_+$) acts over the coordinates in $E_-$ (resp. $E_+$).
\begin{lemma}\label{l:twopot}
We assume ${\rm Dm}(S_-,0),{\rm Dm}(S_+,0)>0$. The Minkowski potential of the divergence-free tensor $S$ given in (\ref{eq:Spm}) is
$$\theta(x_-,x_+)=a_-\theta_-(x_-)+a_+\theta_+(x_+)$$
where
\begin{equation}\label{eq:apm}
a_-=({\rm Dm}(S_-,0))^{-\frac q{d-1}}({\rm Dm}(S_+,0))^{\frac{q-1}{d-1}},\qquad a_+=({\rm Dm}(S_-,0))^{\frac{p-1}{d-1}}({\rm Dm}(S_+,0))^{-\frac p{d-1}}.
\end{equation}
\end{lemma}

\bigskip

\begin{cor}\label{c:Spm} 
The determinantal masses of $S,S_-$ and $S_+$ satisfy the relation
\begin{equation}\label{eq:DmSSS}
{\rm Dm}(S,0)=({\rm Dm}(S_-,0))^{\frac{p-1}{d-1}}({\rm Dm}(S_+,0))^{\frac{q-1}{d-1}}.
\end{equation}
\end{cor}

\bigskip

\bepr

We look for a potential $\theta$ given as a linear combination of $\theta_-$ and $\theta_+$, where we have to identify the coefficients $a_\pm$. 

Because the calculation of determinantal masses requires an approximation procedure (to smooth out the vertex singularity of the potential), we begin by considering smooth convex functions $\xi_\pm$ over $\R^p$ and $\R^q$, instead of $\theta_\pm$. We have easily
$$\widehat{a_-\xi_-+a_+\xi_+}=
\begin{pmatrix} a_-^{p-1}a_+^q(\det{\rm D}_+^2\xi_+)\widehat{{\rm D}_-^2\xi_-} & 0 \\ 0 & a_-^pa_+^{q-1}(\det{\rm D}_-^2\xi_-)\widehat{{\rm D}_+^2\xi_+} \end{pmatrix}.$$
When $\xi_-$ converges uniformly to $\theta_-$, $\det{\rm D}_-^2\xi_-$ tends towards the mesure ${\rm Dm}(S_-;0)\,\delta_{x_-=0}$ while $\widehat{{\rm D}_-^2\xi_-}$ tends to $S_-$. Passing to the limit, we infer
$$\widehat{a_-\theta_-+a_+\theta_+}=
\begin{pmatrix} a_-^{p-1}a_+^q{\rm Dm}(S_+;0)\,S_-\otimes\delta_{x_+=0} & 0 \\ 0 & a_-^pa_+^{q-1}{\rm Dm}(S_-;0)\,\delta_{x_-=0}\otimes S_+ \end{pmatrix}.$$
We recover the tensor $S$ by chosing the solution $(a_-,a_+)$ of the system
$$a_-^{p-1}a_+^q{\rm Dm}(S_+;0)=1,\qquad a_-^pa_+^{q-1}{\rm Dm}(S_-;0)=1.$$
This gives us the formul\ae\, (\ref{eq:apm}).

Finally, the body enclosed by the image of $\nabla(a_-\theta_-+a_+\theta_+)$ is the Cartesian product of the bodies enclosed by the images of $a_\pm\nabla_\pm\theta_\pm$ respectively. In terms of volumes, we have therefore
$${\rm Dm}(S;0)=a_-^pa_+^q{\rm Dm}(S_-,0)\cdot{\rm Dm}(S_+;0),$$
which gives the relation (\ref{eq:DmSSS}).

\enpr

\paragraph{Remarks.}
\begin{itemize}
\item
Formula (\ref{eq:DmSSS}) is a rather natural generalization of the identity $\det M=\det M_-\times\det M_+$ for a block-diagonal matrix.
\item
The procedure above can be generalized to the situation where $\R^d$ is split into an arbitrary number of orthogonal subspaces whose union contains the $\eta_J$'s.
\item The orthogonality between $E_-$ and $E_+$ is not an essential ingredient, because it is always possible to make a linear change of variables which maps isometrically $E_-$ over $\R^p\times\{0\}$ and $E_+$ over $\{0\}\times\R^d$, and to modify $T$ accordingly (see \cite{Ser_IHP}). The general formula (\ref{eq:DmSSS}) will however contain a factor reflecting the angle between $E_-$ and $E_+$.
\end{itemize}

\bigskip

The simplest example of a direct sum is given by the potential
$$\theta_{\rm abs}(x)=\sum_1^d|x_j|=:\|x\|_1.$$
The corresponding tensor is diagonal
$$\Lambda_{\rm abs}=2\sum_1^d\vec e_j\otimes \vec e_j\,\delta_{\R\vec e_j}=2\begin{pmatrix} \delta_{\hat x_1=0} & & \\ & \ddots & \\ & &  \delta_{\hat x_d=0} \end{pmatrix},$$
where as usual $\hat x_j=(x_,\ldots,x_{j-1},x_{j+1},\ldots,x_d)\in\R^{d-1}$. We have immediately
\begin{equation}\label{eq:Mabs}
{\rm Dm}(\Lambda_{\rm abs};0)=2^d,
\end{equation}
which suggests to adopt the convention that
\begin{equation}
\label{eq:convdel}
\left(\delta_{\hat x_1=0} \cdots   \delta_{\hat x_d=0}\right)^{\frac1{d-1}}=\delta_0,
\end{equation}
This extends the well-known formula, when $d=2$, that $\delta_{x_2=0}\delta_{x_1=0}=\delta_0$.

\section{Proof of Theorem \protect\ref{th:main}}\label{s:bin}

The physical domain is $\R^n$. We consider the generic case where the collisions form a discrete set and are only pairwise. We assume without loss of generality that there is no collision at initial time.
We denote $T$ the mass-momentum tensor constructed in Section \ref{s:MMT}. 

To begin with, we choose a time $\tau>0$ at which there is no collision and we set $H_\tau=(0,\tau)\times\R^n$.

\paragraph{A complement to the mass-momentum tensor.} 
Let $K$ be a kink of a trajectory, happening at a point $x^*\in H_\tau$. The incoming/outgoing velocities of the particle under consideration being $v,v'$ respectively, with $v'\ne v$, we complete the free family
$$V=\binom1v,\quad V'=\binom1{v'}$$
into a basis $(V,V',z_2,\ldots,z_n)$ of $\R^{1+n}$. Here $(z_2,\ldots,z_n)$ is some orthonormal basis of ${\rm Span}(V,V')^\bot$. We define the positive semi-definite tensor
$$S_K=\sum_{j=2}^nz_j\otimes z_j\delta_{\sigma_j},$$
where $\sigma_j:=(x^*-\epsilon_K z_j,x^*+\epsilon_K z_j)$ is a segment of direction $z_j$. Mind that the vectors $z_j$  do depend on the kink, even if it is not explicit in our notation. The lengths $\epsilon_K>0$ are small enough that the corresponding segments are contained in $H_\tau$, do not overlap and do not intersect the support of $T$ away from $K$.
Because the number of collisions is finite in $H_\tau$, the entries of $S$ are finite measures.

We form an auxiliary tensor
$$T'=T+S,\qquad S:=\sum_{\hbox{kinks in }H_\tau}b_KS_K,$$
where the positive numbers $b_K$ will be chosen later. It is positive semi-definite, supported by a graph, yet it is not divergence-free, since ${\rm Div}\,T'=\sum b_K{\rm Div}\,S_K\ne0$. Because ${\rm Div}\,S_K$ is a sum of Dirac masses at the end points $x^*\pm \epsilon z_j$, we have instead
$$\|{\rm Div}\,T'\|=2(n-1)\sum_{\hbox{kinks in }H_\tau}b_K.$$

At initial and final time, $T'$ coincides with $T$, and therefore we have
$$\|T'_{0\bullet}(t=0)\|=\sum\sqrt{1+|v(0)|^2\,}$$
where the sum runs over the particles. We infer
$$\|T'_{0\bullet}(t=0)\|\le \sum\left(1+\frac{|v(0)|^2}2\right)\le M+E=N\left(1+\frac12\,\bar v^2\right).$$
Likewise we have $\|T'_{0\bullet}(t=\tau)\|\le M+E$. Since $T'$ is a finite measure, we are therefore in position to apply Compensated Integrability to $T'$ in $H_\tau$.

\bigskip

The left-hand side of (\ref{eq:CIvrtx}) involves the integral of $(\det T')^{\frac1n}$, which vanishes identically, and the determinantal masses at the kinks of the trajectories. A kink $K$ at a point $x^*$ involve the three vectors $V,V'$ and $Q=V-V'$. A combination of Corollaries \ref{c:VWZ}, \ref{c:Spm} and of formula (\ref{eq:Mabs}) yields the following calculation (here we have $p=2$ and $q=n-1$):
$${\rm Dm}(T';x^*)=2^{n-3}b_K^{\frac{n-1}n}|V\wedge V'|^{\frac1n}$$
Writing (\ref{eq:CIvrtx}) for the tensor $T'$, we obtain
$$\sum_{\hbox{kinks in }H_\tau}b_K^{\frac{n-1}n}|V\wedge V'|^{\frac1n}\le c_n\left(2(M+E)+2(n-1)\sum_{\hbox{kinks in }H_\tau}b_K\right)^{1+\frac1n}$$
for another universal constant, still denoted $c_n$. 

We now introduce auxiliary parameters $\lambda>0$ and $\beta_K>0$, and we set $b_K=\lambda\beta_K^{\frac n{n-1}}$. We infer
$$\sum_{\hbox{kinks in }H_\tau}\beta_K|V\wedge V'|^{\frac1n}\le c_n\lambda^{\frac1n-1}\left(M+E+\lambda(n-1)\sum_{\hbox{kinks in }H_\tau}\beta_K^{\frac n{n-1}}\right)^{1+\frac1n}.$$
Choosing
$$\lambda:=(M+E)\left(\sum_{\hbox{kinks in }H_\tau}\beta_K^{\frac n{n-1}}\right)^{-1},$$
we infer
$$\sum_{\hbox{kinks in }H_\tau}\beta_K|V\wedge V'|^{\frac1n}\le c_n(M+E)^{\frac2n}\|\vec\beta\|_{\ell^{\frac n{n-1}}}.$$
The above inequality is valid for every choice of positive parameters $\beta_K$. Since the left-hand side is a scalar product $\langle\vec\beta,\vec D\rangle$, and the dual space of $\ell^{\frac n{n-1}}$ is $\ell^n$, it tells us  that
$$\|\vec D\|_{\ell^n}\le c_n(M+E)^{\frac2n}.$$
In other words, we have
$$\sum_{\hbox{kinks in }H_\tau}|V\wedge V'|\le c_n(M+E)^2.$$
Since
$$|V\wedge V'|=\sqrt{|v'-v|^2+|v\wedge v'|^2\,}\,,$$
we obtain our first estimate
$$\sum_{\hbox{kinks in }H_\tau}\sqrt{|v'-v|^2+|v\wedge v'|^2\,}\le c_n(M+E)^2.$$
Remarking that the right-hand side does not depend upon the time length $\tau$, we actually have
\begin{equation}\label{eq:unbal}
\sum_{\hbox{kinks in }H_\infty}\sqrt{|v'-v|^2+|v\wedge v'|^2\,}\le c_n(M+E)^2,
\end{equation}
where now the sum extends over all the history.

\paragraph{Using the scaling.}
Equation (\ref{eq:unbal}) is not acceptable from a physical point of view. It lacks homogeneity: One should not add a mass and an energy (right-hand side) or two different powers of velocities (left-hand side). To overcome this flaw, we notice that from a given flow, one can construct a one-parameter family of flows, by changing the time scale. This trick was used already in the context of the Euler equations of a compressible fluid, see \cite{Ser_IHP}.

We consider particles of same radius $a$. If $\mu>0$ is given, a trajectory $t\mapsto X(t)$ in the original flow ${\cal F}_1$ gives rise to a trajectory $t\mapsto X_\mu(t):=X(\mu t)$ in the new flow ${\cal F}_\mu$. The velocity is $v_\mu(t)=\mu v(\mu t)$. The flow parameters become
$$M_\mu=M,\qquad E_\mu=\mu^2E,\qquad \bar v_\mu=\mu\bar v.$$
Applying (\ref{eq:unbal}) to ${\cal F}_\mu$ results in a parametrized inequality
$$\sum_{\hbox{kinks in }H_\infty}\sqrt{\mu^2|v'-v|^2+\mu^4|v\wedge v'|^2\,}\le c_n(M+\mu^2E)^2,\qquad\forall\,\mu>0.$$
Choosing $\mu^2=M/2E=\bar v^{\,-2}$, we obtain
$$\sum_{\hbox{kinks in }H_\infty}\sqrt{\bar v\,^2|v'-v|^2+|v\wedge v'|^2\,}\le c_nM^2\bar v\,^2.$$
This is equivalent to (\ref{eq:estbin}) because of 
$$\frac{a+b}2\le\sqrt{a^2+b^2}\le a+b$$
for positive numbers. This ends the proof of Theorem \ref{th:main}.


\begin{thebibliography}{00}


\bibitem{Alex} R. K. Alexander. {\em The infinite hard sphere system.} M.S. thesis, University of California at Berkeley, (1975).


\bibitem{BFK} D. Burago, S. Ferleger, A. Kononenko. Uniform estimates on the number of collisions in semi-dispersing billiards. {\em Annals of Mathemematics,} Second series, {\bf 147} (1998), pp 695--708.


\bibitem{BuIv} D. Burago, S. Ivanov. Examples of exponentially many collisions in a hard ball system.ArXiv preprint, arXiv:1809.02800v1, 2018.


\bibitem{BD} K. Burdzy, M. Duarte. A lower bound for the number of elastic collisions.  {\em Commun. Math. Phys.}, {\bf372} (2019), pp 679--711.


\bibitem{CIP} C. Cercignani, R. Illner, M. Pulvirenti. {\em The mathematical theory of dilute gases.} Springer-Verlag, New York (1994).




\bibitem{Illn} R. Illner. On the number of collisions in a hard sphere particle system in all space. {\em Transport Theory and Stat. Phys.,} {\bf 18} (1989), pp 71--86.


\bibitem{Illde} R. Illner. Finiteness of the number of collisions in a hard sphere particle system in all space II: arbitrary diameters and masses. {\em Transport Theory and Stat. Phys.,} {\bf 19} (1990), pp 573--579.

\bibitem{Min} H. Minkowski. Volumen und Oberfl\"ache. {\em Math. Annalen} {\bf57} (1903), pp 447--495.

\bibitem{Pog} A. V. Pogorelov. {\em The Minkowski multidimensional problem}. Scripta Series in Mathematics. V. H. Winston \& Sons, Washington, D.C.; Halsted Press (John Wiley \& Sons), New York--Toronto--London (1978).




\bibitem{Ser_IHP} D. Serre. Divergence-free positive symmetric tensors  and fluid dynamics. {\em Annales de l'Institut Henri Poincar\'e (analyse non lin\'eaire),} {\bf 35} (2018), pp 1209--1234.

\bibitem{Ser_CI} D. Serre. Compensated integrability.  Applications to the Vlasov--Poisson equation  and other models of mathematical physics. {\em J. Maths. Pures \& Appl.}, {\bf127} (2019), pp 67--88.
 
 
\bibitem{Sin} Ya. Sinai. Hyperbolic billiards. {\em Proceedings of the International Congress of Mathematicians}, Vol. I, II (Kyoto 1990), pp 249--260. Math. Soc. Japan, Tohyo (1991).
 

\bibitem{Vas} L. N. Vaserstein. On systems of particles with finite-range and/or repulsive interaction. {\em Commun. Math. Phys.,} {\bf 69} (1979), pp 31--56.




\end{thebibliography}
\end{document}